\documentclass[12pt]{article}
\usepackage{amssymb}
\usepackage{amsmath}
\usepackage{theorem}

\newcommand{\bC}{{\mathbb C}}
\newcommand{\bQ}{{\mathbb Q}}
\newcommand{\bR}{{\mathbb R}}
\newcommand{\bF}{{\mathbb F}}
\newcommand{\bP}{{\mathbb P}}
\newcommand{\bZ}{{\mathbb Z}}

\newcommand{\cC}{{\mathcal C}}
\newcommand{\cD}{{\mathcal D}}
\newcommand{\cL}{{\mathcal L}}

\newcommand{\cN}{{\mathcal N}}

\newcommand{\cO}{{\mathcal O}}

\newcommand{\cW}{{\mathcal W}}
\newcommand{\ra}{{\rightarrow}}
\newcommand{\amp}{\mathrm{amp}}
\newcommand{\nef}{\mathrm{nef}}
\newcommand{\nod}{\mathrm{nod}}
\newcommand{\Gr}{{\mathrm{Gr}}}
\newcommand{\Pic}{{\mathrm{Pic}}}
\newcommand{\Aut}{{\mathrm{Aut}}}

\newcommand{\Sym}{{\mathrm{Sym}}}
\newcommand{\rk}{{\mathrm{rk}}}
\newcommand{\Span}{{\mathrm{Span}}}

{\theorembodyfont{\rm}
\newtheorem{dfn}{Definition}[section]
\newtheorem{coro}[dfn]{Corollary}
\newtheorem{ques}[dfn]{Question}
\newtheorem{exam}[dfn]{Example}
\newtheorem{rem}[dfn]{Remark}

}

\newtheorem{conj}[dfn]{Conjecture}
\newtheorem{theo}[dfn]{Theorem}
\newtheorem{prop}[dfn]{Proposition}

\author{Brendan Hassett and Yuri Tschinkel}
\title{\Large Rational curves on holomorphic symplectic fourfolds}

\date{March 2001}

\begin{document}

\maketitle

\section{Introduction}

One of the main problems in the theory of 
irreducible holomorphic symplectic manifolds is the description of the ample 
cone in the Picard group.  The goal of this paper is to formulate explicit Hodge-theoretic criteria for the ampleness of line bundles on certain irreducible holomorphic symplectic manifolds.      
It is well known that for K3 surfaces the ample cone is governed by 
$(-2)$-curves.  More generally, we expect that
certain distinguished two-dimensional homology classes 
of the symplectic manifold should correspond to explicit families of rational 
curves, and that these govern its ample cone.  

The program for analyzing the ample cone of a symplectic
manifold divides naturally into three parts.  First, for each deformation 
type of irreducible holomorphic sympectic manifolds we identify distinguished Hodge 
classes in $H_2({\bZ})$ that should be represented by rational curves.  We consider
$H_2({\bZ})$ and $H^2({\bZ})$ 
as quadratic lattices with respect to a natural quadratic form 
(the Beauville form discussed in Section \ref{sect:generalities}) 
and distinguish orbits in $H_2(\bZ)$ under the orthogonal group.  
These orbits are often characterized by the `squares'
of the corresponding elements, i.e., the value of the quadratic form.
These distinguished classes should be in one-to-one correspondence
with certain geometrically described rational curves on $F$.  
In many cases, one can use deformation arguments to show
that, if a given distinguished class 
represents a rational curve of a certain
type then this remains 
true under deformation (see Section \ref{sect:deform}).  
Second, one shows that rational curves with the geometry described govern the ample
cone of $F$.  This entails classifying possible contractions
of symplectic manifolds - a very active topic of current research (see
the discussion of the literature below) - and interpretting this classification
in terms of the numerical properties of the contracted curves.  It also
involves the classification of base loci for sections of line bundles
on a symplectic manifold.  The third part of the program is to show that
a divisor class satisfying certain numerical conditions arises from a big
line bundle and thus yields a birational transformation of the symplectic manifold.  
This aspect of the program is still largely conjectural;  see the work
of Huybrechts cited below.  

We are mainly concerned with the first step of this program in a specific case.     
Let $F$ be an irreducible holomorphic symplectic fourfold
deformation equivalent to the punctual Hilbert scheme $S^{[2]}$ for
some K3 surface $S$. 
Given the Hodge structure on $H^2(F)$, 
we describe explicitly (but conjecturally) the cone of
effective curves on $F$ and, by duality, the ample cone of $F$.  
As in the case of
K3 surfaces, each divisor class of square $-2$ induces
a reflection preserving the Hodge structure.  The `birational
ample cone' is conjectured to be the interior of a fundamental domain 
for this reflection group.  However, the ample cone may be strictly
smaller than the birational ample cone, owing to
the existence of elementary transformations along
$\bP^2$'s in $F$.  The corresponding classes have
square $-10$ with respect to the Beauville form.

\

Here we give a brief and incomplete overview of work on related
problems.  
Wilson has studied Calabi-Yau threefolds from a similar point of view
(see \cite{wilson1},\cite{wilson2},\cite{wilson3}, and \cite{wilson4}.)
If $F$ is a Calabi-Yau threefold then the Picard group of $F$
is equipped with two integer-valued forms: 
a cubic form $\mu$ (the intersection form) and a linear form $c_2(F)$
(obtained by intersecting
with the second Chern class of the tangent bundle.)  
Wilson gives criteria for the existence of birational contractions
and elliptic fibrations in terms of the number-theoretic properties
of these forms.  

Namikawa \cite{namikawa} and Shepherd-Barron \cite{shepherd-barron},
have proven structural results on the geometry of birational morphisms from holomorphic symplectic manifolds.  
There are also results in this direction by Burns, Hu, and Luo
\cite{BHL} and Wierzba \cite{wierzba}.  Matsushita \cite{matsushita1}
\cite{matsushita2} has proven a detailed
description of fiber structures on irreducible holomorphic 
symplectic manifolds.  Huybrechts has conjectured a 
projectivity criterion for irreducible holomorphic symplectic
manifolds and elaborated consequences of his criterion.
See \cite{huybrechts1} and the erratum in \cite{huybrechts2};
we will be careful to distinguish results fully proved in
\cite{huybrechts1} from those which remain conjectural.  
Markman \cite{markman} has developed a theory of generalized elementary
transformations for moduli spaces of sheaves on K3 surfaces.  

\

Our study of rational curves on symplectic manifolds began with
the detailed study of a particular example:
the variety $F$ parametrizing lines on a cubic fourfold $X$ is
an irreducible holomorphic symplectic manifold.  The existence of rational
curves on $F$ coincides with the presence of rational ruled
surfaces on $X$.  Our conjectures therefore shed light on
the effectivity of certain codimension-two cycles on $X$.
Conversely, the projective geometry of
cubic fourfolds provides a useful laboratory where
we may test our claims.  As an application of our conjectures,  
we find that the presence of distinguished Hodge classes 
on $X$ often implies the existence of special unirational
parametrizations of $X$.

\

This paper is organized as follows.  In Section \ref{sect:generalities}
we recall basic results and conjectures for irreducible
holomorphic symplectic manifolds.  In the next section, we introduce
the notion of {\em nodal classes} and state 
our main conjectures.  In Section \ref{sect:deform} we give
some deformation-theoretic evidence for our conjectures.  
The rest of the paper is devoted to examples supporting 
the conjectures. Section \ref{sect:symm} is devoted to
Hilbert schemes $S^{[2]}$ for K3 surfaces of small degree.  
We describe examples of nonnodal rational curves and certain
codimension-two behavior in Section \ref{sect:nonnodal}.
We turn to the projective geometry of cubic fourfolds
in the last section.  Questions of rationality and
unirationality are addressed in Section \ref{sect:unirat}.

Throughout, we work over ${\bC}$.  A {\em primitive} element
of an abelian group $A$ is one that cannot be written in the
form $nx$ for any $x\in A$ or $n\in {\bZ},n>1$.  An {\em
indecomposable}
element of a monoid is an element which cannot be written in the form
$a+b$ for some nonzero $a$ and $b$ in this monoid. 
Recall that an element $v$ of a convex real cone $C$ is an
{\em extremal ray} if, for any $u,w \in C$ with $u+w=v$ we 
necessarily have $u,w \in {\bR}_+v$.

\

{\bf Acknowledgements.} The first author
was partially supported by an NSF postdoctoral fellowship,
NSF continuing research grant 0070537, and
The Institute of Mathematical
Sciences at the Chinese University of Hong Kong.  
The second author was partially supported by the NSA. 
We thank D. Matsushita and Y. Namikawa for sending us their preprints. 
Our treatment of rationality questions in Section
\ref{sect:unirat} benefitted from discussions 
with D. Saltman and I. Dolgachev.

\section{Generalities}
\label{sect:generalities}

Let $F$ be an irreducible, holomorphic symplectic manifold of dimension $2n$.  
This means that $F$ is compact, K\"ahler, simply connected, and 
$H^0(F, \Omega^2_F)$ is spanned by an everywhere-nondegenerate 2-form
(cf. \cite{beauville-82} and \cite{bogomolov}).  
The second cohomology group $H^2(F,\bZ)$  
carries a nondegenerate integer-valued 
quadratic form $\left(,\right)$,
the {\em Beauville form}.  It has signature $(3,b_2(F)-3)$ and
its restriction to $H^{1,1} \cap H^2(F,{\bZ})$ has signature
$(1,b_2(F)-3)$,  
where $b_2(F)$ is the second Betti number
(see \cite{beauville-83} and
\cite{huybrechts1} \S 1.9 for more details).  

Using the universal coefficient theorem, we can
extend the Beauville form to a ${\bQ}$-valued form on
$H_2(F,{\bZ})$.  Concretely,
given some primitive $R\in H_2(F,{\bZ})$,
there exists a unique class $w\in H^2(F,{\bQ})$ such that
$
Rv=\left(w,v\right)
$
for all $v\in H^2(F,\bZ)$.  We set $\left(R,R\right)=\left(w,w\right)$.
Let $\rho \in H^2(F,{\bZ})$ denote
the primitive class such that $c\rho=w$ for some $c>0$.
Note that $R$ is of type $(2n-1,2n-1)$ iff $\rho$ is of type
$(1,1)$.  Conversely, given a primitive $\rho \in H^2(F,{\bZ})$ 
with $\left(\rho,H^2(F,{\bZ}) \right)=d{\bZ}$ and $d>0$, there exists a
primitive class $R\in H_2(F,{\bZ})$ with $dRv=\left(v,\rho \right)$ for
all $v\in H^2(F,{\bZ})$.  

Throughout this paper, the square of a divisor class means the
square with respect to the Beauville form.
In the sequel we will assume that $F$ has a polarization $g$;
note that $\left(g,g\right)>0$ \cite{huybrechts1} \S 1.9.  
Denote by 
$$
\Pic_+(F,g) = \{v\in \Pic(F)\,|\, \left(v,g\right) >0\}
$$
the positive halfspace (with respect to $g$ and the Beauville form).  
Let $\Lambda_+(F,g)\subset \Pic_+(F,g)$ be the vectors with positive 
square.  

We denote by $\Lambda_{\amp}(F)$ and $\Lambda_{\nef}(F)$
the monoids of ample and nef divisor classes. 
Let $N_1(F)$ be the group of classes of 1-cycles (up to numerical 
equivalence), $NE(F)\subset N_1(F)_{\bR}$
the cone of effective curves, and $\overline{NE}(F)$ its
closure.  
We denote by $\Lambda_+^*(F,g)$ the set of classes $R$
such that the corresponding $\rho$ is contained in $\Lambda_+(F,g)$.

We next review properties of line bundles on 
polarized irreducible holomorphic symplectic manifolds which
follow from standard results 
of the minimal model program.  

\begin{prop}\label{prop:standard}
Let $(F,g)$ be a polarized irreducible holomorphic symplectic manifold.
\begin{enumerate}

\item 
A class $\lambda \in \Pic(F)$ is ample iff $\lambda \in \Lambda_+(F,g)$
and $\lambda C>0$
for each curve $C\subset F$ (see \cite{huybrechts1} Cor. 6.4).

\item 
Any class $\lambda \in \Lambda_{\nef}(F)$ which is big
has the property that the line bundle $\cL(\lambda)$
has no cohomology and $\cL(m\lambda)$ is globally generated for $m\gg 0$;
it therefore defines
a birational morphism $b:F \ra Y$
(see \cite{KMM}, Remark 3-1-2 and Theorem 1-2-3).

\end{enumerate}
\end{prop}
The following statements were conjectured by Huybrechts \cite{huybrechts2}
(but stated as Theorem 3.11 and Corollary 3.10 in \cite{huybrechts1}
with an incomplete proof):

\begin{conj}
\label{conj:huybrechts}
Let $F$ be an irreducible holomorphic symplectic manifold.  
\begin{enumerate}
\item
$F$ is projective
iff there exists a class $g \in \Pic(F)$ with 
$\left(g,g\right)>0$.

\item
If $g$ is a polarization for $F$ then
any class $\lambda\in \Lambda_+(F,g)$ is big.
\end{enumerate}
\end{conj}

Now assume $F=S^{[n]}$, the Hilbert scheme 
of length $n$ subschemes of a K3 surface $S$.
Then we have an isomorphism 
$$\Pic(S^{[n]}) \simeq \Pic(S) \oplus_{\perp} {\bZ}e$$
compatible with the Beauville form (see \cite{beauville-83}).  
Each divisor $f$ on $S$ determines a divisor on $S^{[n]}$,
also denoted by $f$, and
corresponding to the subschemes with some support in $f$.
The locus of subschemes with support at 
fewer than $n$ points has class $2e$ and $\left(e,e\right)=-2(n-1)$.

More generally, if $F$ is 
deformation equivalent to $S^{[n]}$ then
the Beauville form on $L:=H^2(F,{\bZ})$ is an even,
integral form isomorphic to 
$$
U^{\oplus 3}  
\oplus_{\perp} (-E_8)^{\oplus 2} \oplus_{\perp} {\bZ}e,
$$
where $U$ is a hyperbolic plane, 
$E_8$ the positive-definite quadratic
form associated to the corresponding 
Dynkin diagram.

\begin{prop}
Assume that $n=2$ so that 
$$L\simeq U^{\oplus 3}  
\oplus_{\perp} (-E_8)^{\oplus 2} \oplus_{\perp} (-2).$$
The orbits of primitive elements $v\in L$ under the 
action of $\Gamma=\Aut(L)$ are classified by 
$\left( v,v\right)$ and the ideal $\left( v,L \right)$
(which equals $\bZ$ or $2\bZ$.)
\end{prop}
{\em Proof.}
We shall classify primitive imbeddings of the lattice $K:=\bZ v=(2d)$
(where $\left( v,v \right)=2d$) into $L$.  For simplicity, 
we first restrict to the case $d\neq 0$.  The basic technical
tool is the {\em discriminant group} $d(L):=L^*/L$ and the
associated $\bQ/2\bZ$-valued quadratic form $q_L$ \cite{nikulin}.
We have $d(L)\simeq \bZ/2\bZ$ with $q_L$ equal to $-\frac{1}{2}\pmod{2\bZ}$
on the generator.
Let $K^{\perp}$ denote
the orthogonal complement to $K$ in $L$, $d(K)$ and $d(K^{\perp})$
the discriminant groups, and $q_K$ and $q_{K^{\perp}}$ the
corresponding forms, so that
$d(K)\simeq \bZ/2d\bZ$ with $q_K$ equal to
$\frac{1}{2d}\pmod{2\bZ}$ on the generator $\frac{v}{2d}$.
We have the sequence of inclusions
$$K\oplus K^{\perp} \subset L \subset L^* \subset K^* \oplus (K^{\perp})^*;$$
$L^*$ consists of the elements of $K^*\oplus (K^{\perp})^*$
which are integral on $L$.  Let $H$ (resp. $H^*$)
be the image of $L$ (resp. $L^*$) in $d(K)\oplus d(K^*)$, so
that $d(L)=H^*/H$.  Note that $H$ is isotropic and $H^*$
is the annihilator of $H$ with respect to
$q_K\oplus q_{K^{\perp}}$ (or the
${\bQ}/{\bZ}$-valued bilinear form associated to it).
The projection of $H$ into $d(K)$ is
injective and its image is a subgroup of index one or two, depending
on whether $\left(v,L\right)=\bZ$ or $2\bZ$.  

In the first case, $d(K^{\perp})$ contains the projection of $H$ as an index 
two subgroup.  Since $H$ is isotropic for $q_K\oplus q_{K^{\perp}}$,
the restriction $q_{{K^{\perp}}}|H$ takes value
$-\frac{1}{2d}\pmod{2\bZ}$ on the generator.  
Let $x\in H^*$ be a nontrivial element projecting
to $0$ in $d(K)$, which may be regarded as an element
of $d(K^{\perp})$.  We have $q_{K^{\perp}}(x)=-\frac{1}{2}\pmod{2\bZ}$ because
$x$ generates $d(L)$.  
Since $x$ annihilates $H$, it follows that
$d(K^{\perp})=H+\bZ x\simeq \bZ /2d\bZ\oplus_{\perp} \bZ/2\bZ$.  This 
determines the discriminant form (and signature) of $K^{\perp}$ 
completely, which determines it up to isomorphism
\cite{nikulin} 1.14.2.  The classification of the primitive imbeddings 
$K\hookrightarrow L$ follows from \cite{nikulin} 1.15.1.

The proof in the second case is similar but easier,
as $d(K^{\perp})$ is equal to $H$, so its discriminant form
is easily computed.  

If $d=0$, we may produce an element $e'\in L$ with $\left(v,e'\right)=0$,
$\left(e',e'\right)=-2,$ and $\left(e',L\right)=2\bZ$.  By our
previous argument
$e'\in \Gamma e$, so we may assume 
that $w\in U^{\oplus 3}\oplus (-E_8)^{\oplus 2}$.
Then the result follows from \cite{lo-peters} \S 2. \hfill $\square$

\

When $n=2$, we have an identity 
$c_2(F)\cdot v \cdot v=30\left(v,v \right)$ (given
in \cite{huybrechts1} 1.11 up to a multiplicative constant).  
Hence Riemann-Roch takes the form
$$\chi(F,\cL(v))=\frac{1}{8}(\left(v,v\right)+4)(\left(v,v\right) +6).$$

\section{Conjectures}
\label{sect:conj}

In this section $(F,g)$ is a $g$-polarized 
irreducible holomorphic symplectic manifold, deformation equivalent to the
Hilbert scheme of length-two subschemes of a K3 surface.  
Let $E$ be the (possibly infinite!) set of classes
$\rho \in \Pic_+(F,g)$
satisfying one of the following:
\begin{enumerate}
\item $\left(\rho,\rho  \right)=-2$ and $\left(\rho,L\right)=2\bZ$,
\item $\left(\rho,\rho  \right)=-2$ and $\left(\rho,L\right)=\bZ$,
\item $\left(\rho,\rho  \right)=-10$ and $\left(\rho,L\right)=2\bZ$,
\end{enumerate} 
Let $E^*$ be the corresponding classes $R\in H_2(F,{\bZ})$;  this 
means that for some $\rho \in E$ we have
$$\left( v,\rho \right)=\begin{cases} 
			Rv & \text{ where } \left(\rho,L\right)={\bZ} \\
		       2Rv & \text{ where } \left(\rho,L \right)=2{\bZ}
		       \end{cases}
$$
for each $v\in L$.  In particular, $R$ satisfies one of the following
\begin{enumerate}
\item $\left(R,R  \right)=-\frac{1}{2}$, 
\item $\left(R,R  \right)=-2$, 
\item $\left(R,R  \right)=-\frac{5}{2}$.
\end{enumerate} 

Let $N_E(F,g)\subset H_2(F,\bZ)$ be the 
smallest real cone containing $E^*$ and the elements  $R\in N_1(F)$
such that $R\cdot g>0$ and the corresonding $\rho$ has nonnegative
square. 
Note that 
the boundary of $N_E(F,g)$ is polyhedral in a neighborhood of any
boundary point with negative square.  
This follows from the fact that the Beauville form is negative definite
on the orthogonal complement to $g$ in $\Pic(F)$.

We now can state our main conjecture.  
\begin{conj}[Effective curves conjecture]
\label{conj:main}
$$NE(F)=N_E(F,g).$$
\end{conj}
The analogous theorem for K3 surfaces may be found
in \S 1.6 and 1.7 of \cite{lo-peters} 
(see also \cite{beauville-85}).

The classes in $E^*$ that are extremal in (the closure of) $N_E(F,g)$
will be called {\em nodal classes}, 
in analogy with the 
terminology for K3 surfaces (see \cite{lo-peters},
Section 1.4).  The nodal classes are denoted $E^*_{\nod}$ and the
corresponding classes in $E$ are denoted $E_{\nod}$.  
Since the boundary of $N_E(F,g)$ is polyhedral in a
neighborhood of any nodal class, it follows that a
class $R\in E^*$ is nodal iff no positive multiple
of $R$ is decomposable.

\begin{rem} Consider the monoid $N_E(F,g)\cap H_2(F,\bZ)$.
In contrast to the situation for K3 surfaces, it is possible 
for a $(-2)$-class in the monoid to be indecomposable
but not nodal.  Indeed, there exist examples with
$\rk \, N_1(F)=2$.  
\end{rem}

Conjecture \ref{conj:main} and Proposition
\ref{prop:standard} yield a characterization of
the ample monoid:
\begin{eqnarray*}
\Lambda_{\amp}(F)&=&\{ \lambda\in \Lambda_+(F,g):
\left(\lambda,v \right) > 0 \text{ for each }\\
& & v\in \Pic_+(F,g)
\text{ with } \left(v,v\right)\ge 0 \text{ or } v\in E\}.
\end{eqnarray*}
The signature of the Beauville form implies that if 
$\left(\lambda,\lambda\right)$ and $\left( v,v \right)$ are both
nonnegative then $\left(\lambda,v \right)\ge 0$;
indeed if $\lambda $ and $v$ are linearly independent then 
strict inequality follows.   
Furthermore, verifying the positivity condition against nonextremal
classes is clearly redundant.  We therefore obtain the following
simplification:
\begin{conj}[Consequence of Conjecture \ref{conj:main}]
\label{conj:ample}
\begin{eqnarray*}
\Lambda_{\amp}(F)&=&\{ \lambda\in \Lambda_+(F,g):
\left(\lambda,v \right) > 0 \text{ for each }
v\in E_{\nod}\}.
\end{eqnarray*}
\end{conj}
This generalizes Proposition 1.9 of \cite{lo-peters}.

We digress to consider the monoid of effective divisors. Its 
description is the same as for K3 surfaces.

\begin{conj}[Effective divisors conjecture]
The monoid of effective divisors is generated by the
elements of
$\Pic_+(F,g)$ with square $\ge -2$.  
\end{conj}

\begin{rem}
A $(-2)$-class is extremal in the
(conjectured) cone of effective divisors 
iff it is indecomposable in the (conjectured) monoid of 
effective divisors.  We expect that each of these
$(-2)$-classes is realized by a conic bundle over a K3 surface.  
\end{rem}

We next discuss classes $\lambda\in \Lambda_{\nef}(F)$,
i.e., those in the boundary of the ample cone.  
Conjecture \ref{conj:huybrechts} implies that
when $\left(\lambda,\lambda \right)>0$ 
the sections of $\cL(m\lambda)$
for $m\gg 0$ give a birational morphism $b:F \ra Y$.  
Any curve represented by a nodal class
$R\in E^*_{\nod}$ orthogonal to $\lambda$ is contracted
by $b$.

\begin{conj}[Nodal classes conjecture]\label{conj:nodal}
Each nodal class $R\in E^*_{\nod}$ represents a rational curve contracted
by a birational morphism $b$ given by sections of $\cL(m\lambda)$,
where $\lambda$ is any class on the boundary of the
ample cone satisfying $\lambda R=0$.  
\begin{enumerate}
\item If $\left(R,R\right)=-\frac{1}{2},-2$
(i.e., the corresponding $\rho$ is a $(-2)$-class) then $\rho$ 
is represented by a family of rational curves parametrized
by a K3 surface.
This family can be blown down to rational double points. 
\item If $\left(R,R\right)=-\frac{5}{2}$
(i.e., the corresponding $\rho$ is a $(-10)$-class) then $\rho$ 
is represented by a family of lines
contained in a plane $\bP^2$. 
This plane can be contracted to a point.
\end{enumerate} 
\end{conj}

The following theorem of Namikawa provides support for this
conjecture (part of the theorem was
proved first by Shepherd-Barron \cite{shepherd-barron} and related results
were obtained by Wierzba \cite{wierzba}).

\begin{theo}(\cite{namikawa} Props. 1.1 and 1.4, \cite{shepherd-barron})
\label{thm:namikawa}
Let $b\, :\, F\ra Y$ be a birational projective morphism from a 
projective holomorphic
symplectic manifold to a normal variety. 
\begin{enumerate}
\item There exists a subvariety $Z\subset Y$ of codimension at least
four,
so that $Y\setminus Z$ is
singular along a smooth codimension-two subvariety $S$, which
admits a nondegenerate holomorphic 
two-form.  Furthermore, $Y\setminus Z$ has rational double points of
fixed type along each connected component of $S$. 
\item Assume that $P\in Y$ is an isolated 
Gorenstein singular point and
that at least one irreducible component of 
$G:=b^{-1}(P)$ is normal.
Then $G\simeq \bP^n$ with normal bundle $\Omega^1_{\bP^n}$.    
\end{enumerate}
\end{theo}

A key tool in the study of K3 surfaces is the {\em Weyl group},
the group generated by reflections with respect to the $(-2)$-classes
in the Picard group (see \cite{beauville-85} or \cite{lo-peters}).  
Our conjectures imply that there is an analog in higher dimensions.  
For each $\rho\in E$ with square $-2$, we obtain
a reflection $s_{\rho}$ given by the formula
$$s_{\rho}(v)=v+\left(v,\rho \right)\rho.$$  
The group $\cW$ generated by these reflections is called the 
{\em generalized Weyl group}.  
Let $C_+(F,g)$ be the smallest real cone containing $\Lambda_+(F,g)$
and $\cD(F,g)$ the subcone of $C_+(F,g)$
defined by $\left(v,\rho\right)\ge 0$ for each
$\rho \in E$ with square $-2$.  This is a fundamental domain for
the action of $\cW$ on $C_+(F,g)$.  The $(-10)$-classes
in $E$ are walls for a subdivision of $\cD(F,g)$ into subchambers.  
The interior of each
subchamber is the ample cone for a
symplectic birational model for $F$;
the subchamber containing $g$ is the ample cone of $F$.  
In other words, $\cD(F,g)$ is the closure of the birational
ample cone of $F$.

One does not generally
expect that the closure of the ample cone should be a fundamental domain
for the action of a group $\cW' \supsetneq \cW$.  Algebraically,
$(-10)$-classes do not yield reflections.  
Geometrically, one does not expect to find a group relating
the various birational models of $F$.  For a concrete example,
see the discussion below of the Fano variety of lines 
on a cubic fourfold of discriminant eight. 

\

We now describe the square-zero classes on the boundary
of the closure of the ample cone:

\begin{conj}\label{conj:abelian}
Let $\lambda$ be a primitive square-zero class on the boundary
of the closure of the ample cone.  
Then the corresponding line bundle
$\cL=\cL(\lambda)$ has no higher cohomology and its
sections yield a morphism
$$a\, :\, F \ra {\bP}^2$$
whose generic fiber is an abelian surface.  
\end{conj}

The following theorem of Matsushita describes fiber space structures
on holomorphic symplectic manifolds. 

\begin{theo}(\cite{matsushita1})\label{theo:matsushita} 
Let $F$ be a projective irreducible holomorphic symplectic manifold and 
$a\, :\, F\ra B$ be a fiber space structure with
normal base $B$ and generic fiber $F_b$.  Then we have:
\begin{enumerate}
\item $a$ is equidimensional in codimension-two points of $B$.
\item $K_{F_b}$ is trivial and there 
exists an abelian variety $\tilde{F}_b$ 
and an \'etale morphism $\tilde{F}_b\ra F_b$.
\item $\dim(B)=n$ and $B$ has  $\bQ$-factorial log-terminal
singularities.
\item  $-K_B$ is ample and $B$ has Picard number 1.
\item The polarization on $B$ pulls back to a 
square-zero divisor on $F$. 
\end{enumerate}   
Furthermore, if $F$ has dimension four then $F_b$ is
an abelian surface. 
\end{theo}
Matsushita has also proved that such fibrations are Lagrangian
\cite{matsushita2}.  For convenience, we provide a proof in our
special case:

\begin{prop}\label{prop:lagran}
In addition to the hypotheses of Theorem \ref{theo:matsushita} assume
that the dimension of $F$ is four. Then the fibers $F_b$ are Lagrangian. 
\end{prop}

{\em Proof.} 
It suffices to prove this result for smooth fibers. 
Recall the formula for the Beauville form (see \cite{huybrechts1}, 1.9)
$$
(\alpha, \alpha)=\int(\sigma\overline{\sigma})\alpha^2
-\left(\int\sigma\overline{\sigma}^{2}\alpha\right)
\left(\int\sigma^2\overline{\sigma}\alpha\right),
$$
where $\sigma$ is the generator of $H^0(F,\Omega^2_F)$,
normalized so that $\int (\sigma {\overline \sigma})^2=1$.  
If $\alpha $ is the pullback of the polarization on $B$ then the class
of $F_b$ is equal to some rational multiple $c\alpha^2$. 
Since $\alpha$ is of type $(1,1)$, type considerations imply that
the second term vanishes. On the other hand, the first term is equal to
$$
\frac{1}{c}\int_{F_b}\sigma|_{F_b}\overline{\sigma}|_{F_b}. 
$$
Since $(\alpha, \alpha)$ 
is zero by the previous theorem it follows that $\sigma$ 
restricted to $F_b$ is also zero and therefore $F_b$ is Lagrangian
(here we use the Hodge-Riemann bilinear relations).\hfill $\square$

\begin{rem}
It is instructive to compare Conjecture \ref{conj:abelian} 
with Wilson's results.  Let $F$ be a Calabi-Yau
threefold and $D$ a nef divisor class corresponding
to a nonsingular point of the cubic hypersurface $\mu_F=0$ and satisfying
$Dc_2(F)\ne 0$. Then for some $n>0$ the linear series $|nD|$ is free 
and induces an elliptic fiber space structure on $F$ 
(see \cite{wilson1} \S 3 and \cite{wilson2} \S 1).  The assumption
$c_2(F)D\ne 0$ can be weakened under further technical hypotheses (see
\cite{wilson3} and \cite{wilson4}).  
\end{rem}

\begin{rem}
Let $F$ be an irreducible holomorphic symplectic
manifold and $\lambda\in \Lambda_{\nef}(F)$ a divisor
with $\left( \lambda, \lambda \right)=0$.  Is $\lambda$ necessarily
semiample?  
\end{rem}

\section{Deformation theory}
\label{sect:deform}

\subsection{Deformations of subvarieties}
\label{subsect:ran}

In this section, we work with arbitrary irreducible holomorphic symplectic
manifolds. 

\begin{theo}\label{theo:def}
Let $F$ be an irreducible
holomorphic symplectic manifold of dimension $2n$ and $Y$ a
submanifold of dimension $k$. 
Assume either that $Y$ is Lagrangian, or that all
of the following hold:
$
\cN_{Y/F}= \Omega^1_Y\oplus \cO^{\oplus 2n-2k}_Y,
$
the restriction of the symplectic form to $Y$ is zero,
and $H^1(\cO_Y)=0$.   
Then the deformation space of $Y$ in $F$ is smooth
of dimension $2n-2k$. 

The deformations of $F$ arising from deformations of 
the pair $(Y,F)$ are precisely those preserving the
sub-Hodge structure 
$$
{\rm ker}(H^{*,*}(F)\rightarrow H^{*,*}(Y)).
$$
\end{theo}

{\em Proof.}
This is very similar to the proof of 
Corollary 3.4 in \cite{ran} (see \cite{voisin2} for
the Lagrangian case). 
By \cite{ran}, Corollary 3.2, any obstructions to deforming
$Y$ in $F$ lie in the kernel of the natural projection 
$$
\pi_{0,2}\, : \, H^1(N_{Y/F}) \rightarrow H^1(\Omega^1_Y) \otimes
\mathrm{ker}(H^0(\Omega^2_F) \rightarrow H^0(\Omega^2_Y))^*.
$$
Under our hypotheses this map is an isomorphism.    \hfill
$\square$

\begin{coro}\label{coro:def}
Keep the hypotheses of Theorem \ref{theo:def}. 
Assume furthermore that the cohomology of $Y$ is 
generated by divisor classes. Then the 
deformations of $F$ arising from a deformation of
the pair $(Y,F)$ are precisely those for which the 
image of 
$$
H_2(Y,\bZ)\rightarrow H_2(F,\bZ)
$$
remains algebraic. 
\end{coro}

\subsection{Applications} 

In this section we will assume that $F$ is deformation
equivalent to $S^{[2]}$ for a K3 surface $S$.
We show that the locus where our
conjectures hold is open (in the analytic topology) in the
moduli space.  We will find 
points in the moduli space where our conjectures
hold in subsequent sections.

\begin{theo}
Let $F$ be as above and $R$ a nodal class on $F$.
Assume that Conjecture \ref{conj:nodal} holds for $R$ and
$F$.  Let $F'$ be a small projective deformation of $F$
such that $R$ deforms to a class $R'$ of type $(3,3)$.  
Then the Conjecture remains true for $F'$ and $R'$.  
\end{theo}

{\em Proof.} 
Let $\pi \,:\, {\cal F}\ra \Delta$ be a deformation of $F$ over a disc and 
${\cal A}$ a nef and big line bundle on ${\cal F}$
such that ${\cal A}$ has degree zero on $R$.
By Corollary \ref{coro:def}, $R$ deforms to a family of
rational curves ${\cal  R}\subset {\cal F}$ over $\Delta$. 
In each fiber  ${\cal R}_t$ itself deforms in a two parameter family.  
This family is parametrized
by a deformation ${\cal S}_t$  of the K3 surface $S$ 
if $\rho$ is a $(-2)$-class.
It is parametrized by a $\bP^2$ if $\rho$ is a $(-10)$-class.
Clearly, these are contracted by the sections of some power of ${\cal A}$.
\hfill $\square$

\begin{theo}
Let $F$ be as above and $\lambda$ be a nef, square-zero
divisor class on $F$.
Assume that Conjecture \ref{conj:abelian} holds for $\lambda$ and
$F$.  Let $F'$ be a small projective deformation of $F$
such that $\lambda$ deforms to a divisor class $\lambda'$.
Then the Conjecture remains true for $F'$ and $\lambda'$.  
\end{theo}

{\em Proof.}
This is a consequence of Theorem \ref{theo:def} and 
the Lagrangian property proved in \ref{prop:lagran}.
\hfill $\square$

\subsection{Examples}

In this section we give examples of submanifolds satisfying the
conditions of Theorem \ref{theo:def}.  We also refer the
reader to \cite{voisin2} for examples in the Lagrangian case.   
Assume that $Y$ is a complete homogeneous space
under a reductive algebraic
group or a toric variety, and assume that 
the normal bundle to $Y$ is 
of the form stated above. 
Theorem \ref{theo:def} 
shows that the deformation space of $Y$ in $F$
is smooth.  Moreover, the locus in the deformation
space of $F$ corresponding to manifolds containing
a deformation of $Y$ has codimension 
equal to the rank
of the N\'eron-Severi group of $Y$.  
Here are some specific examples:

\begin{exam}
Given $Y={\bP}^n\subset F$, the deformations
of $F$ containing ${\bP}^n$ form a divisor in the deformation space.
For instance, if $S$ is a K3 surface containing a smooth rational
curve $C$ then $C^{[n]}\simeq \bP^n \subset S^{[n]}$. 
Let $R\in H_2(F,{\bZ})$ be the class of a line in $\bP^n$ 
and $\rho \in H^2(F,{\bZ})$ the corresponding divisor 
(i.e. $2Rv=\left(\rho,v\right)$ for $v\in H^2(F,{\bZ})$).  
Then $\rho=2C-e$ and $\left(\rho,\rho \right)=-2(n+3)$.  
\end{exam}

\begin{exam}
Again, let $S$ be a K3 surface containing a smooth rational curve
$C$ and $F=S^{[n]}$.  Consider the subschemes of $S$ of 
length $n$ with some support along $C$.  The
generic such subscheme is
the union of a point of $C$ and a length $n-1$ subscheme
disjoint from $C$.  Thus we get a divisor $D_1\subset S^{[n]}$
birational to a $\bP^1$ bundle over $S^{[n-1]}$. 
The normal bundle to a generic fiber $Y$ of this bundle
is $\Omega^1_Y\oplus \cO_Y^{\oplus 2n-2}$.  Let
$R=[Y]\in H_2(F,{\bZ})$ and $\rho=[D_1] \in H^2(F,{\bZ})$;
note that $R$ corresponds
to $\rho$ and $\left(\rho,\rho \right)=(C\cdot C)_S=-2$.
The deformations of $S^{[n]}$ containing a deformation 
of $D_1$ are those for which $R$ (or $\rho$)
remains algebraic.  They have
codimension one in the deformation space.  
\end{exam}

The relevance of these examples to Conjecture \ref{conj:nodal}
is discussed in Section \ref{subsect:sing}.  
We will give further examples
in Section \ref{sect:nonnodal}, where we consider 
cases where $F$ is deformation equivalent to $S^{[2]}$
and $Y={\bF}_0,{\bF}_1,$ or ${\bF}_4.$
We digress to give one further example that is particularly
interesting:
\begin{exam} We give a geometric realization of 
certain Mukai isogenies between K3 surfaces (cf. \cite{mukai}).  
Let $S_8$ be a generic degree 8 K3 surface.  
In particular, we assume $S_8$ is realized
as a complete intersection of 3 quadric hypersurfaces in $\bP^5$, 
and the discriminant curve for these quadrics is a smooth 
sextic plane curve $B$.  It follows that each such quadric $Q$ has rank
five or six, and the corresponding family of maximal isotropic
subspaces in $Q$ is parametrized by ${\bP}^3$ or a 
disjoint union two copies of ${\bP}^3$ respectively.  As we
vary $Q$, the families of 
maximal isotropic subspaces are parametrized by
a K3 surface $S_2$ of degree 2, the double cover of ${\bP}^2$
branched over $B$.  
Thus we obtain an \'etale $\bP^3$-bundle ${\cal E}\ra S_2$, 
mapping into $S_8^{[4]}$, with fibers 
satisfying the conditions of Proposition \ref{theo:def}.

This yields an elegant universal construction of Brauer-Severi
varieties representing certain 2-torsion elements of the 
Brauer group of a degree two K3 surface.  Other 2-torsion 
elements are realized as \'etale $\bP^1$-bundles 
${\cal E} \ra S_2$ arising from families of nodal 
rational curves (see the discussion of cubic fourfolds of
discriminant $8$ in Examples \ref{exam:disc8} and \ref{exam:T20}).  
The relationship between Mukai isogenies and Brauer groups
is explored more systematically in the upcoming
thesis of Caldararu \cite{caldararu}.  
\end{exam}

\section{Symmetric squares of K3 surfaces}
\label{sect:symm}

Let $S_{2n}$ be a K3 surface with Picard group generated by
a polarization $f_{2n}$ of degree $2n$.  
The Beauville form restricted to the Picard group takes the form
$$
\begin{array}{c|rr}
 & f_{2n}& e \\
 \hline 
f_{2n} & 2n & 0 \\
e & 0 & -2
\end{array} 
$$
The effective divisor with class $2e$ is called the diagonal.
It is isomorphic to a $\bP^1$-bundle over $S_{2n}$;
the fibers are nodal rational curves.  It follows that an
ample line bundle has class of the form $xf_{2n}-ye$ with $x,y>0$
and Proposition \ref{prop:standard} implies that $2nx^2-2y^2>0$.   
The conjectures in Section \ref{sect:conj} 
give sufficient conditions on $x$ and $y$ for
$xf_{2n}-ye$ to be ample.     

\begin{prop}
Assume that $S_{2n}$ is a K3 surface with a polarization $f_{2n}$
which  embeds 
$S_{2n}$ as a subvariety of $\bP^{n+1}$. 
The line bundle $af_{2n}-e$ on $S_{2n}^{[2]}$ 
is ample whenever $a>1$ or $a=1$ and $S_{2n}$ does not contain 
a line. In particular, $f_{2n}$ lies on the boundary of the 
closure of the ample cone. 
\end{prop}

{\em Proof.}
Let $S$ be a smooth surface embedded in projective space $\bP^r$ 
and not containing a line. Then there is a morphism 
from  the Hilbert scheme $S^{[2]}$ to the Grassmannian $\Gr(2,r)$.
This morphism is finite onto its image. Therefore, the pullback 
of the polarization on the Grassmannian to $S^{[2]}$  is ample.
We apply this to the image of $S_{2n}$ under the line bundle $af_{2n}$.
\hfill $\square$
 
\begin{rem}
In the event that $S_{2n}$ does contain a line $\ell \subset \bP^{n+1}$
the line bundle $f_{2n}-e$ fails to be ample. However, it is nef
and big and a sufficiently high multiple of it gives
a birational morphism contracting the plane $\ell^{[2]}$ and 
inducing an isomorphism on the complement to this plane. In particular, 
there is a nodal $(-10)$-class $2[\ell] -e$ orthogonal to $f_{2n}-e$.
\end{rem}

\subsection{Degree $2$ K3 surfaces}

A K3 surface $S_2$ of degree two can be realized as a double cover
of $\bP^2$ ramified in a curve of degree 6. The quadratic
form $2x^2-2y^2$ does represent $-2$ and $-10$. 
The corresponding nodal classes are $e$ and $2f_2-3e$. 
The second class corresponds to the plane in $S_2^{[2]}$
arising from the double cover. 
Our conjectures predict that the ample cone consists
of classes $xf_2-ye$ where $x,y>0$ and $2x-3y>0$.
The quadratic form also represents 0, but the corresponding 
class $f_2-e$ satisfies  $(2f_2-3e,f_2-e)=-2$. After flopping the
plane the proper transform of $f_2-e$ {\em does} 
yield an abelian surface
fibration (the Jacobian fibration) - as expected. 

Here is a sketch proof
that the class $3f_2-2e$ is nef and big.
Indeed, $3f_2$ is very ample and embeds $S_2$
into $\bP^{10}$. The image is cut out by  quadrics ${\cal I}(2)$. 
Each pair of points on $S_2$ determines a line $\ell$.
The quadrics vanishing on that line form a hyperplane in ${\cal I}(2)$. 
This induces a morphism from $S_2^{[2]}$ to $\bP^{27}=\bP({\cal I}(2)^*)$
given by the sections of the line bundle.

\subsection{Degree $4$ K3 surfaces}

Let $S_4$ be a K3 surface 
with Picard group generated by a polarization of degree 4.
We take $f_4-e$ as the polarization of $S_4^{[2]}$. 
Now we describe the $(-2)$ and $(-10)$-lattice vectors 
in $\bZ f_{4}\oplus \bZ e$ and 
determine which are nodal classes. 
In fact, there are no $(-10)$-lattice vectors.  
The $(-2)$-vectors are of the form $\pm a_m f_4 \mp b_m e$, 
where $a_m\sqrt{2}+b_m=(2\sqrt{2}+3)^m$. The vectors
in the positive halfspace $\Pic_+(S^{[2]}_4,f_4-e)$ 
satisfy $2x-y>0$.  
The nodal classes are  $2f_4-3e$ and $e$. 
It is easy to see that all the other $(-2)$-classes
in the positive halfspace are
decomposable.  We therefore predict that the ample cone
is the interior of the cone spanned by $f_4$ and
$3f_4-4e$. Indeed, $S_4^{[2]}$ has an involution exchanging
$f_4$ and $3f_4-4e$ (given $p,q\in S_4$ the line spanned
by $p$ and $q$ meets $S_4$ in two more points). 

\subsection{Degree $8$ K3 surfaces}

Let $S_8$ be a K3 surface 
with Picard group generated by a polarization of degree 8.
This is the smallest degree case where there are no nodal classes
besides the diagonal.  
Indeed, the quadratic form $8x^2-2y^2$ represents $-2$ and $-10$
exactly when $(x,y)=(0,\pm 1)$ and $(\pm 1,\pm 3)$.  
However, the parity condition for nodal classes of
square $-10$ is not satisfied by $f_8-3e$, i.e., 
$\left(f_8-3e,H^2(F,\bZ) \right) \ne 2\bZ$.  
Therefore, our conjectures imply that the ample cone is the
interior of the cone spanned by $f_8$ and $f_8-2e$,
and the second line bundle yields an abelian surface fibration
$a:F\ra {\bP}^2.$
We have already seen that the ample cone 
is contained in this cone.
For an explicit construction of the abelian surface
fibration, see \cite{hassett-tschinkel} \S 7.  There it
is shown that the symmetric square of a generic K3 surface 
of degree $2n^2$ ($n>1$) admits an abelian surface fibration.

\begin{rem}
This is a counterexample to the theorem in Section 2, p. 463 of 
\cite{markushevich-1}. There it is claimed that 
$S^{[2]}$ of a K3 surface $S$ admits a (Lagrangian) abelian surface
fibration if and only if $S$ is elliptic. 
\end{rem}

\subsection{K3 surfaces containing a rational curve}
\label{subsect:sing}

Let $S$ be a K3 surface containing a rational curve $C$
and let $T$ be the surface obtained by blowing down $C$.
Of course, $T$ has one rational double point.  Consider the map
$b\,:\, S^{[2]}\ra \Sym^2(T)$. 
This map contracts rational curves corresponding to both
$(-2)$ and $(-10)$-nodal classes. 

The Hilbert scheme $S^{[2]}$ contains
a plane $C^{[2]}$ and 
two distinguished divisors $D_1$ and $D_2$, 
birational to $\bP^1$-bundles
over $S$.  The divisor $D_1$ is the locus of
subschemes with some support in $C$ and $D_2$ is the diagonal.
The map $b$ contracts $D_1$ and $D_2$ to surfaces isomorphic
to $T$ and $C^{[2]}$ to the point $p$ where these
surfaces intersect.  The fiber $b^{-1}(p)$ is the union of 
$C^{[2]}$ and $\bF_4$
(cf. Theorem \ref{thm:namikawa}). 

The divisors $D_1$ and $D_2$ have classes $C$ and
$2e$ respectively.  If $R$ is the class of a line in 
$C^{[2]}\simeq {\bP}^2$
then the corresponding divisor class $\rho=2C-e$.
The $(-2)$-class $e$ and the $(-10)$-class $\rho$ are nodal;
the class $C$ is not nodal.  Of course, it becomes nodal
upon flopping $C^{[2]}\subset S^{[2]}$.

\section{Nonnodal smooth rational curves}
\label{sect:nonnodal}

It is well known that for K3 surfaces all smooth rational curves
are nodal and correspond to indecomposable (-2)-classes. 
In Section \ref{subsect:sing} we gave examples
of nonnodal smooth rational curves; these curves were parametrized
by a K3 surface.  
Here we discuss further examples of nonnodal 
smooth rational curves.
As we shall see, these
curves need not be parametrized
by a K3 surface or a $\bP^2$.

We first consider three examples where smooth rational curves do not
correspond to nodal classes, but still correspond to classes
with negative square. 
Let $F=S^{[2]}$, where $S$ is a K3 surface which is a double cover
of a rational  surface $\Sigma$ with Picard group of rank 2.
Then $F$ contains a surface isomorphic to $\Sigma$. 
We emphasize that the results of Subsection \ref{subsect:ran} apply in
this case.  This suggests certain refinements to
Conjecture \ref{conj:nodal}, which we formulate in each example.  
 
\begin{exam}\label{exam:F0}
Let $S\ra \Sigma=\bF_0$ be branched over 
a general curve of type $(4,4)$.  Hence the rulings induce
two elliptic fibrations $E_1$ and $E_2$ which generate the Picard
group and intersect as follows:
$$
\begin{array}{c|rr}
 & E_1 & E_2 \\
 \hline 
E_1 & 0 & 2 \\
E_2 & 2 & 0
\end{array}.
$$
Let $R_1$ and $R_2$ denote the rulings of $\Sigma\subset S^{[2]}$,
with $\rho_1$ and $\rho_2$ their Poincar\'e duals in $\Pic(S^{[2]})$.
We have $\rho_1=E_1 -e $ and $\rho_2=E_2-e$ so that the Beauville form
may be written
$$
\begin{array}{c|rr}
 & \rho_1 & \rho_2 \\
 \hline 
\rho_1 & -2 & 0 \\
\rho_2 & 0 & -2
\end{array}. 
$$
Moreover, the $\rho_i$ generate a saturated sublattice
of the Picard group and 
$$
(\rho_i,H^2(F,\bZ))=\bZ.
$$

The smooth curves in the class $R_1+R_2$ move 
in a 3-parameter family
on $\Sigma\subset F$. However, $\rho_1+\rho_2$ is 
not a nodal class. 
We conjecture that any holomorphic symplectic fourfold deformation
equivalent to a symmetric square of a K3 surface 
with 2 nodal classes $\rho_1$ and $\rho_2$ as above
should contain a surface $\Sigma=\bF_0$. 
\end{exam}

\begin{exam}
Let $S\ra \Sigma=\bF_1$ be branched over a general curve of type
$6R_0+4R_{-1}$, where $R_0$ is the class of the ruling 
and $R_{-1}$ is the class of the exceptional curve. 
The ruling induces an elliptic fibration $E$ on $S$ and the exceptional
curve yields a rational curve $C\subset S$;  these generate
the Picard group and intersect as follows:
$$
\begin{array}{c|rr}
 & E & C \\
 \hline 
E & 0 & 2 \\
C & 2 & -2
\end{array}.
$$

Let $\rho_0$ (resp. $\rho_{-1}$) be the Poincar\'e dual to 
$R_0$ (resp. $2R_{-1}$).  We have $\rho_0=E-e$ and $\rho_{-1}=2C-e$,
so that $\rho_0$ and $\rho_{-1}$ generate a saturated sublattice
on which the Beauville form may be written
$$
\begin{array}{c|rr}
 & \rho_0 & \rho_{-1} \\
 \hline 
\rho_0 & -2 & 2 \\
\rho_{-1} & 2 & -10
\end{array}. 
$$
Moreover, $(\rho_0,H^2(F,\bZ))=\bZ$ and 
$(\rho_{-1},H^2(F,\bZ))=2\bZ$. 

The smooth curves in the 
class $2R_0+R_{-1}$ move in a 4-parameter family
on $\Sigma\subset F$. 
However, $4\rho_0+\rho_{-1}$ is not a nodal class. 
In this case we conjecture that any
$F$ whose cohomology contains  
2 nodal classes $\rho_0$ and $\rho_{-1}$ as above
should contain a surface $\Sigma=\bF_1$. Furthermore, 
we expect that $F$ is a specialization of a variety containing 
a plane $\Pi$ which corresponds to
a $(-10)$-nodal class. This class is equal to 
$2\rho_0+\rho_{-1}$ and 
$\Pi$ specializes to a union of a 
$\bP^2$ and the $\bF_1$ in $F$. 
\end{exam}

\begin{exam}
Let $S\ra \Sigma=\bF_4$ be branched over the union of 
a general curve of type $12R_0+3R_{-4}$ and $R_{-4}$, 
where $R_0$ is the class of the ruling 
and $R_{-4}$ is the class of the exceptional curve. 
Again, the Picard group of $S$ is generated by an elliptic
fibration $S$ and a rational curve $C$ which intersect as follows
$$
\begin{array}{c|rr}
 & E & C \\
 \hline 
E & 0 & 1 \\
C & 1 & -2
\end{array}.
$$
Let $\rho_0$ (resp. $\rho_{-4}$) be the Poincar\'e dual to 
$R_0$ (resp. $R_{-4}$).  Then we have $\rho_0=E-e$ and 
$\rho_{-4}=2C+e$, so $\rho_0$ and $\rho_{-4}$ generate a
saturated sublattice with Beauville form
$$
\begin{array}{c|rr}
 & \rho_0 & \rho_{-4} \\
 \hline 
\rho_0 & -2 & 4 \\
\rho_{-4} & 4 & -10
\end{array}. 
$$
Moreover, $(\rho_0,H^2(F,\bZ))=\bZ$ and 
$(\rho_{-4},H^2(F,\bZ))=2\bZ$. 

The smooth curves in the class 
$5R_0+R_{-4}$ move in a 7-parameter family
on $\Sigma\subset F$. However, 
$5\rho_0+\rho_{-4}$ is not a nodal class. 
In this case we conjecture that any
$F$ with cohomology  
containing  2 nodal classes $\rho_0$ and $\rho_{-4}$ as above
should contain a surface $\Sigma=\bF_4$. We also
expect that $F$ is a specialization of a variety containing 
a plane $\Pi$ which corresponds to
a $(-10)$-nodal class. This class is equal to $4\rho_0+ \rho_{-4}$ and 
$\Pi$ specializes to a union of a $\bP^2$ and the $\bF_4$ in $F$. 
\end{exam}

Next we consider examples of smooth rational curves in $F$
where the corresponding class $\rho$ is of positive
square. 

\begin{exam}
Let $S_2$ be a general  K3 surface of degree 2 
with polarization $f_2$. 
Let $C\subset S_2$ be a rational curve
with two ordinary double points contained in the linear series
$|f_2|$. Let $F\ra \bP^2$ be the compactified Jacobian 
for $|f_2|$. The fibers corresponding to $C$ are isomorphic to
a product of nodal curves with normalization
$\bP^1\times \bP^1$. 
Smooth curves of type $(1,1)$ in $\bP^1\times \bP^1$ 
yield smooth rational curves on $F$, deforming in a 3-parameter
family. The homology class of these rational curves is 
double the class of the curve of type $(1,0)$, and is therefore
not primitive. 
\end{exam}
\begin{exam}
Let $S_4\subset \bP^3$ be a general K3 surface of degree 
4 and $f_4$ its polarization. 
Let $C\subset S_4$ be an elliptic curve in $|f_4|$ with 
two ordinary double points.  Note that 
$C^{[2]}$ is a nonnormal  ruled surface. Its fibers
are smooth rational curves such that the corresponding class
$\rho$ has square 2. Thus we get smooth rational curves
in primitive homology classes such that
the corresponding class $\rho$ has positive square as well. 
\end{exam}

As the rank of the Picard group of $F$ increases we expect
more and more examples of nonnodal smooth rational curves
parametrized by varieties of dimension $>2$. 

\begin{rem}
Let $R\subset F$ be a smooth rational curve with primitive
homology class. 
Then the Hilbert scheme of 
flat deformations of $R$
need not be irreducible and
may have arbitrarily large dimension.  
Take $R\subset \bF_0\subset F$ of bidegree $(1,n)$. 
\end{rem}

\begin{ques}
Assume that $\rk\, \Pic(F)=1$. 
Does there exist a smooth rational curve  on 
$F$? Can we take its class to be  primitive?  
\end{ques}

\section{Cubic fourfolds}
\label{sect:cubics}

In this section, a cubic fourfold generally denotes a smooth cubic
hypersurface $X\subset {\bP}^5$.  
The variety $F$ parametrizing lines on $X$ 
is sometimes called the `Fano variety of lines' - not to be confused
with a variety with ample anticanonical class.  It is known that $F$
is an irreducible holomorphic symplectic fourfold 
deformation equivalent to the Hilbert scheme of length-two
subschemes of a K3 surface
\cite{AK} \cite{beauville-donagi}.  Consequently, the conjectures
of Section \ref{sect:conj} apply.  The existence
of smooth rational curves $R\subset F$ translates into 
the existence of scrolls $T\subset X$.  By definition, a scroll
is the union of the lines parametrized by a smooth rational
curve in the Grassmannian;  it may
have singularities.    
Our conjectures yield simple and 
verifiable predictions for the existence and nonexistence
of scrolls in various homology classes of $X$. 
The presence of these
scrolls yields unirational parametrizations of $X$ of various degrees.

\subsection{Lattices, Nodal Curves, and Scrolls}
\label{sect:genlat}

We recall standard facts about cubic fourfolds.
We say that a cubic fourfold is {\em special} if it contains an
algebraic surface not homologous to any multiple of the square
of the hyperplane class $h^2$.  
Note that the intersection form $\left<,\right>$ on
the primitive cohomology takes the form
$$
(h^2)^{\perp}\simeq \left(\begin{array}{cc} 2& 1\\ 1 &2\end{array}
\right) \oplus_{\perp} U^{\oplus 2}\oplus_{\perp}E_8^{\oplus 2}
$$
(see \cite{hassett-96}, \cite{hassett-96-0}, \cite{beauville-donagi}).
Let $K={\bZ}h^2+{\bZ}T$ be a
saturated sublattice of algebraic classes in the middle cohomology
of $X$.   Then the {\em discriminant} 
$d=d(X,K)$ is the discriminant of $K$.
It is a positive
integer, congruent to $0$ or $2$ modulo $6$.
The special cubic fourfolds
of discriminant $d$ form an irreducible divisor $\cC_d$ in the moduli
space $\cC$ of cubic fourfolds;  $\cC_d$ is nonempty iff $d>6$.  For
instance, $\cC_8$ corresponds to the cubic fourfolds containing
a plane $T_1$ and $\cC_{14}$ 
corresponds to the cubic fourfolds containing
a smooth quartic scroll $T_4$.

The cohomology of a cubic fourfold and 
its Fano variety are closely related
(see \cite{beauville-donagi} for most of what follows).  
The incidence correspondence between 
$X$ and $F$ induces the Abel-Jacobi
map 
$$
\alpha \, :\, H^4(X,\bZ) \ra H^2(F,\bZ),
$$
respecting the Hodge structures.  
We have that 
$\left(\alpha(h^2),\alpha(h^2)\right)=2\left< h^2,h^2\right>$
and 
$$\left(\alpha(v),\alpha(w)\right)=-\left<v,w\right>$$
for $v,w$ primitive.  Note that $g:=\alpha(h^2)$
is the polarization on $F$ induced from the Grassmannian.
The incidence correspondence induces a second map
$$\beta\,:\,H^6(F,\bZ)\ra H^4(X,{\bZ})$$
respecting the Hodge structures.  We can compose
to obtain
$$\psi:H_2(F,\bZ) \ra H^6(F,\bZ) \stackrel{\beta}{\ra} H^4(X,\bZ) 
\stackrel{\alpha}{\ra} H^2(F,\bZ) \ra H_2(F,\bZ),$$
where the first map is Poincar\'e duality and the last
map is induced by the Beauville form.  We
have $\psi(g)=2g$ and $\psi(v)=-v$
for $v$ orthogonal to $g$.  

Suppose that $F$ contains a smooth rational curve $R$
of degree $n$.  Let $\tilde T$ be the universal line
restricted to $R$ and $T\subset X$ the corresponding scroll
sweeped out by $R$, which also has degree $n$.  Note that
the formula
$\left<T,\Sigma \right>=R\cdot \alpha(\Sigma)$
(for $\Sigma \in H^4(X,\bZ)$)
follows from the incidence correspondence.
Combining this with our computation of $\psi$, we obtain
$$\left<T,T \right>=R\cdot \alpha(T)=\left(R,\psi(R) \right)=
\frac{n^2}{2}-\left(R,R\right).$$

We use $T_{n,\Delta}$ to denote a scroll $T$ of degree $n$ 
for which the map
$$\tilde T \ra T\subset X$$
has singularities equivalent 
to $\Delta$ ordinary double points
(by definition, $\Delta$ is the number given by the double point formula). 
A Chern class computation gives
\begin{equation}
\left< T_{n,\Delta},T_{n,\Delta} \right>=3n-2+2\Delta \label{eqn:selfint}
\end{equation}
and we obtain the formula
\begin{equation}
\Delta=\frac{1}{4}(n^2-6n+4-2\left(R,R\right)). \label{eqn:DeltaR}
\end{equation}

The lattice generated by $h^2$ and $T_{n,\Delta}$ has discriminant
\begin{eqnarray*}
d(n,\Delta)&=&3(3n-2+2\Delta)-n^2=6\Delta-(n^2-9n+6) \\
&=& \frac{n^2}{2}-3\left(R,R\right).
\end{eqnarray*}
This lattice has discriminant $>6$,
so we obtain the lower bound
$$
\Delta\ge \Delta_{\min}(n):=\lceil \frac{1}{6}(n^2-9n+6)+1 \rceil.
$$
In particular, a cubic fourfold cannot contain smooth scrolls
of degree $>7$.  

\begin{rem}
The lattice ${\bZ}h^2+{\bZ}T_{n,\Delta}$ need not be saturated.  
For instance, if $n=8$ and
$\Delta=5$ then $d(8,5)=32$.  However, the lattice generated
by $h^2$ and $T_{8,5}$ has index $2$ in its saturation.  
\end{rem}

\begin{prop}\label{prop:delta}
Let $X$ be a cubic fourfold, with Fano variety $F$.
Let $R\subset F$ be a nodal rational curve and
$T_{n,\Delta}$ the corresponding scroll.
Then $\Delta$ takes the following values:
$$\Delta =\begin{cases}
(m-2)(m-1) \text{ if $n=2m$}; \\
(m-1)^2   \,\, \text{and}\,\, m(m-2)  
\text{ if $n=2m+1$}.  
\end{cases}
$$
\end{prop}

{\em Proof.} This is a consequence of Equations 
\ref{eqn:selfint} and \ref{eqn:DeltaR} above.  We observe that
$n$ is even when $\left(R,R\right)=-2$ and $n$ is odd
when $\left(R,R\right)=-\frac{1}{2}$ or $-\frac{5}{2}$.  
\hfill $\square$ 

We summarize the numerical predictions for {\em nodal} 
scrolls of small degree in the following table:
$$
\begin{array}{c|rrrrrrrrrrrrrr}
n      & 2 & 3 & 4 & 5 & 5 & 6 & 7 & 7 & 8 & 9 & 9 & 10 & 11 & 11 \\
\hline
\Delta & 0 & 0 & 0 & 0 & 1 & 2 & 3 & 4 & 6 & 8 & 9 & 12 & 15 & 16 \\
d(n,\Delta)      & 8 &12 &14 &14 & 20& 24& 26& 32& 38& 42& 48& 56 & 62 & 68 
\end{array}
$$
\begin{rem} 
We can obtain cubic fourfolds containing scrolls $T_{n,\Delta}$
with more double points by exploiting 
nonnodal smooth rational curves on the corresponding Fano
variety (see Examples \ref{exam:T41} and \ref{exam:T63}).
\end{rem}

\subsection{Unirational parametrizations}
\label{sect:unirat}
We start with a classical example:  if $X$ is a cubic
fourfold containing a smooth quartic scroll $T_{4,0}$
then $X$ is rational.  One would like
to generalize this construction to other special cubic
fourfolds.  

\begin{prop} \label{prop:unirat}
Let $X$ be a cubic fourfold with Fano variety $F$.
Assume that $F$ contains a smooth rational curve $R$
of degree $n$, with corresponding scroll $T_{n,\Delta}$.  
Assume that this corresponding scroll $T$ is not a cone.
Then there exists a rational map
$$\phi:\bP^4 \dashrightarrow X$$
with 
$$
\deg(\phi)
=\binom{n-2}{2}-\Delta 
=\frac{(n-2)^2}{4}+\frac{\left(R,R \right)}{2}+1.
$$
\end{prop}

{\em Proof.}  Our assumptions imply that
$R$ parametrizes pairwise disjoint 
lines in $X$.  Given generic $\ell_1,\ell_2$, the cubic
surface 
$$\Span(\ell_1,\ell_2)\cap X$$
contains two disjoint lines and thus is rational.
We therefore obtain a cubic surface bundle
$$\begin{array}{ccc}
Y           & \stackrel{\sigma}{\rightarrow} &\Sym^2(R)\simeq{\bP}^2 \\
\stackrel{\phi}{}\downarrow \  & & \\
X & & 
\end{array}$$
so that the fiber over the generic point contains two disjoint lines.
Consequently, $Y$ is rational over $\bP^2$ and thus is
a rational variety.

To compute the degree of $\phi$, it suffices to compute
the number of double points arising from a generic
projection of the scroll $T$ into $\bP^4$.  The map
${\tilde T} \ra \bP^4$ has singularities equivalent
to $\binom{n-2}{2}$ double points;  $\Delta$ of
these are from the singularities of $T$.  
We obtain the second formula for $\deg(\phi)$ by applying
the Equation \ref{eqn:DeltaR} of Section \ref{sect:genlat}.  
\hfill $\square$

\

This demonstrates that the existence of rational curves
in certain homology classes of $F$ implies that $X$ is rational. 
Unfortunately, our conjectures indicate that such rational
curves are quite rare.  If $R$ is nodal then $\left(R,R \right)\ge
-\frac{5}{2}$, so $\deg(\phi)=1$ only when $n=4$ (see also Examples 
\ref{exam:disc26} and \ref{exam:T52}).

However, we do obtain some interesting new unirational
parametrizations of cubic fourfolds.  Recall that in
\cite{hassett-96}, the Fano variety of lines 
on the generic cubic fourfold of discriminant 
$2(N^2+N+1)$ ($N>1$) was shown to be isomorphic to $S^{[2]}$ of a
K3 surface $S$.  In particular, it contains
nodal rational curves $R$ of degree $2N+1$ 
with $\left(R,R\right)=-\frac{1}{2}$.  One can show that the 
scroll corresponding to a generic such curve is not a cone, 
hence Proposition \ref{prop:unirat} applies.  We obtain
$$\deg(\phi)=N^2-N+1,$$
which is always odd.  
In particular, the cubic fourfolds with {\em odd} degree
unirational parametrizations are dense in the moduli space.

Cubic fourfolds are known to admit unirational parametrizations
of degree two.  Thus the cubic fourfolds described above admit
unirational parametrizations of relatively prime degrees.  
There are few examples of irrational varieties with this property.
Many common invariants used to detect irrationality (like the
unramified cohomology of the function field) vanish 
in this situation.

\subsection{Cubic fourfolds of small discriminant}
\label{sect:smalldisc}

In this section we specialize our conjectures
to Fano varieties of lines on {\em general} special cubic
fourfolds of discriminant $d$. We obtain predictions
on the existence and nonexistence of 
scrolls $T_{n,\Delta}$ on $X_d\in \cC_d$.   
We verify these predictions in Section \ref{sect:data}. 
Throughout we write $g=\alpha(h^2)$ and $\tau = \alpha(T)$.
\begin{exam}[$d=8$]\label{exam:disc8}
For $X_8\in \cC_8$ (resp. $F_8$) we have
intersection pairing (resp. Beauville form): 
$$
\begin{array}{c|rr}
 & h^2& T \\
          \hline 
h^2 & 3 & 1 \\
T   & 1 & 3
\end{array} \,\, \hskip 1cm
\begin{array}{c|rr}
 & g & \tau \\
 \hline 
g      & 6 & 2 \\
\tau   & 2 & -2
\end{array},
$$
so $\tau $ is a $(-2)$-class (note that
$(\tau, H^2(F_8,\bZ))=\bZ$.) 
There is also a $(-10)$-class: $\rho=g-2\tau$.  One can check that
these classes are nodal.  Therefore our conjectures
predict a plane in $\Pi \subset F_8$ whose lines have degree one
in the Grassmannian.  This corresponds to a plane in $X_8$. 
They also predict a family of 
rational curves in $F_8$ parametrized by 
a K3 surface which correspond to quadric cones in $X_8$
(see Example \ref{exam:T20}).

This example illustrates our previous discussion
concerning the action of the Weyl group.  Here we have
$$C_+(F_8,g)=\{ ag-b\tau: 3a+b>0,  \ a-b>0 \}$$
and the fundamental domain for the action of the Weyl group is
$$\cD(F_8,g)=\{ ag-b\tau: a+b\ge 0, \  a-b\ge 0 \}.$$
The conjectures predict that the ample cone should be 
$$\Lambda_{\amp}(F_8)=\{  ag-b\tau: a+b>0, \  a-3b > 0 \};$$
the nef cone is bounded by clases of square $0$ and $64$.
If $F'$ denotes the elementary transformation of $F_8$ along 
the plane $\Pi$, we expect
$$\Lambda_{\amp}(F')=\{  ag-b\tau:  a-b>0, \  -a+3b>0 \};$$
the nef cone is bounded by classes of square $8$ and $64$.  
In particular, the two subchambers of $\cD(F_8,g)$ are not
conjugate.  
\end{exam}

\begin{exam}[$d=12$]\label{exam:disc12}
For  $X_{12}$ and  $F_{12}$ we have
pairings: 
$$
\begin{array}{c|rr}
 & h^2& T \\
          \hline 
h^2 & 3 & 3 \\
T   & 3 & 7
\end{array} \,\,\hskip 1cm  
\begin{array}{c|rr}
 & g & \tau \\
 \hline 
g      & 6 & 6 \\
\tau   & 6 & 2
\end{array}.
$$
The $(-10)$-classes are given by $2\tau-g $ and $3g-2\tau$. 
Our conjectures predict that $F_{12}$ contains two 
projective planes. The lines on these planes correspond to 
families of cubic scrolls on $X_{12}$ (see Example \ref{exam:T30}).
\end{exam}
\begin{exam}[$d=14$]\label{exam:disc14}
For  $X_{14}$ and  $F_{14}$ we have
pairings: 
$$
\begin{array}{c|rr}
 & h^2& T \\
          \hline 
h^2 & 3 & 4 \\
T   & 4 & 10
\end{array} \,\, \hskip 1cm
\begin{array}{c|rr}
 & g & \tau \\
 \hline 
g      & 6 & 8\\
\tau   & 8 & 6
\end{array}.
$$
The nodal classes classes are given by
$2g-\tau$ and $2\tau-g$. 
Note that $(2g-\tau, H^2(F_{14},\bZ))=\bZ$ and
$(2\tau-g, H^2(F_{14},\bZ))=2\bZ$.
The first corresponds to a family of rational curves
of degree 4 on $F_{14}$ parametrized by a K3 surface.
The second corresponds to a family of rational curves
of degree 5 also parametrized by a K3 surface
(see Examples \ref{exam:T40} and \ref{exam:T50}).
\end{exam}
\begin{exam}[$d=20$]\label{exam:disc20}
For $X_{20}$ and $F_{20}$ we have pairings
$$
\begin{array}{c|rr}
 & h^2 & V \\
\hline 
h^2 & 3 & 4 \\
V & 4 & 12
\end{array} \, \,\hskip 1cm 
\begin{array}{c|rr}
 & g & v \\
\hline 
g & 6 & 8 \\
v & 8 & 4 
\end{array},  
$$
where $v=\alpha(V)$. 
There are no $(-2)$-classes but there are two nodal $(-10)$-classes:
$e_1=2v-g, e_2=19g-8v$. The corresponding rational curves
on $F_{20}$ have
degrees 5 and 25, respectively
(cf. Example \ref{exam:T51}).
There is an involution interchanging $e_1$ and $e_2$ given by: 
$$
\begin{array}{ccr}
g    & \mapsto & 5g-2 v \\
v & \mapsto & 12g-5v
\end{array}.
$$
\end{exam}
\begin{exam}[$d=26$]\label{exam:disc26}
For $X_{26}$ and $F_{26}$ we have pairings
$$
\begin{array}{c|rr}
 & h^2 & T \\
\hline 
h^2 & 3 & 5 \\
T & 5 & 17
\end{array} \, \, \hskip 1cm
\begin{array}{c|rr}
 & g & \tau \\
\hline 
g & 6 & 10 \\
\tau & 10 & 8 
\end{array}.  
$$
This lattice does not represent $-10$.  The nodal 
$(-2)$-classes are $2\tau -g$ and $109g-38\tau$.  Note
that 
$$(2\tau-g,H^2(F_{26},\bZ))=(109g-38\tau,
H^2(F_{26},\bZ))=2\bZ.$$  
Our conjecture predicts two families of rational curves 
parametrized by K3 surfaces, with degrees 
$7$ and $137$ respectively.  

We next apply our conjecture on effective classes 
to derive the nonexistence of a $T_{5,2}$ on $X_{26}$.
By Proposition \ref{prop:unirat}, the existence of such a surface would 
imply the rationality of $X_{26}$. 
Let us assume that $T_{5,2}\subset X_{26}$  with ruling
$R$. We may take $T$ for the class of $T_{5,2}$. We compute
the class $\rho$ corresponding to $R$. Since 
$$
\frac{1}{2}\rho\cdot g = 5 \,\,\hskip 1cm 
\frac{1}{2}\rho\cdot \tau =17
$$
we get $\rho=5g-2\tau$. If we write 
$$
\rho= a (2\tau -g) + b(109g-38\tau)
$$
then we find that $a = -7/45$ and $b=2/45$.  This implies
that $R$ is not contained in the (conjectured) monoid 
of effective classes.  
We shall show in Example \ref{exam:T52} that 
a quintic scroll in $\bP^5$ cannot have two double points.
\end{exam}

\begin{ques}
How can one systematize the argument for the nonexistence
of $T_{5,2}$'s on a (general) $X_{26}$? 
More precisely, let 
$X$ be a cubic fourfold containing a scroll $T_{n,\Delta}$ and 
assume that the lattice containing $h^2$ and 
$T_{n,\Delta}$ generates the lattice of algebraic classes in 
$H^4(X,\bZ)$. Do the values 
obtained in Proposition \ref{prop:delta} give
upper bounds for $\Delta$ in terms of $n$? 
\end{ques}

\subsection{Data}
\label{sect:data}

In this section we present data from projective
geometry concerning the existence of 
scrolls on cubic fourfolds. 
We organize the information by the degree of the
scroll. 

First of all, let us observe that
a scroll of degree $n$ with ordinary double points can be obtained
by projecting a smooth nondegenerate scroll of degree $n$ in
$\bP^{n+1}$  from a suitable linear subspace. 

\begin{exam}[$T_{2,0}$] \label{exam:T20}
Observe that a scroll of degree two  cannot have
ordinary double points at all.  It is easy to see that
the general cubic fourfold of discriminant 
8 contains such a scroll. 
Furthermore, these scrolls are parametrized 
by a K3 surface of degree
2 (see  \cite{voisin}, \cite{hassett-96}, \cite{hassett-99}).  
\end{exam}
\begin{exam}[$T_{3,0}$] \label{exam:T30}
A scroll of degree 3 also does not have
any ordinary  double points and it is 
contained in a general cubic
fourfold of discriminant 12. 
On a fixed cubic fourfold these scrolls 
are parametrized by {\em two} disjoint 
$\bP^2$'s;  given one scroll $T$, there is a residual scroll
$T'$ obtained by intersecting a linear and a quadratic
hypersurface containing $T$ (see \cite{hassett-96-0}). 
These correspond to two distinct $(-10)$-classes. 
\end{exam}
\begin{exam}[$T_{4,0}$] \label{exam:T40}
A nondegenerate scroll of degree 4 in $\bP^5$
does not have any ordinary double points. A general cubic
fourfold of discriminant 14 contains  a family 
of such scrolls, parametrized by a smooth 
K3 surface of degree 14. 
The corresponding class is a nodal $(-2)$-class
(see, for example, \cite{beauville-donagi},
\cite{hassett-96-0}). 
\end{exam}
\begin{exam}[$T_{4,1}$] \label{exam:T41}
This example is closely related to Example \ref{exam:F0}.
We will explain why the locus
of cubic fourfolds containing 
a quartic scroll with one ordinary double point
has codimension 2 in moduli. 
Consider a cubic
fourfold $X$ containing such a scroll $T_{4,1}$.
Note that $T_{4,1}$ is degenerate and is contained in a  
singular cubic threefold
$Y$. We specialize first to the case where the quartic scroll
degenerates to the union of two quadric scrolls.  Each of these
quadric scrolls is residual to a plane, and these planes intersect
at a single point.  
What can we say about $Y$ in this case?  A cubic
threefold containing two such planes is obtained as follows.
Let $C$ be a genus 4 stable curve obtained by taking a curve $C_1$
of type $(1,3)$ on a quadric 
surface $Q$, along with the union of two rulings
$C_2$ and $C_3$ of type $(1,0)$.  Note that $C$ is canonically
imbedded in $\bP^3$.  Then $Y$ is the image of $\bP^3$
under the linear series $|L|$
of cubics cutting out $C$.  
The planes in $Y$ are the total 
transforms of the lines $C_2$ and $C_3$.
We claim that $Y$ contains a family of scrolls $T_{4,1}$,
parametrized by $\bP^3$.
In particular, a cubic fourfold containing two planes meeting at a
point also contains a three parameter family of $T_{4,1}$'s.  
These are obtained by taking the proper transforms of the quadric
surfaces $Z$ in $\bP^3$ containing 
the lines $C_2$ and $C_3$.  These form
a linear series with projective dimension three.  The restriction
of $|L|$ to $Z$ is a linear series 
of type $(1,3)$ with two base points
(i.e., the points of $Z \cap C_1$ not lying on $C_2$ or $C_3$).  
The image of $Z$ is a $T_{4,1}$.

This corresponds to a situation where $F$ contains a 
surface isomorphic to $\bP^1 \times \bP^1$.  The 
hyperplane sections give a 3-parameter family of rational
curves on $F$.  As we have seen, such fourfolds should 
lie in codimension two (see \ref{theo:def}). 
\end{exam}
\begin{exam}[$T_{5,0}$] \label{exam:T50}
The cubic fourfolds of discriminant 14 
also contain a family of quintic scrolls, 
parametrized by {\em the same} K3 surface which para\-met\-rizes
the quartic scrolls.  
The corresponding class
is a second nodal $(-2)$-class
(see \cite{beauville-donagi} or
\cite{hassett-96-0}). 
\end{exam}
\begin{exam}[$T_{5,1}$] \label{exam:T51}
A general cubic fourfold $X_{20}$ of discriminant 20 
contains a family of $T_{5,1}$'s parametrized by a $\bP^2$. 
It is known that $X_{20}$ contains a Veronese surface $V$.
This also follows from Theorem \ref{theo:def} once we obtain
a $\bP^2\subset F_{20}$.  
The conic curves in $X_{20}$ lying in
$V$ are parametrized by $\bP^2$ as well.  For each such
curve $C$, let $H$ be the plane spanned by $C$ so that
$$X_{20} \cap H=C \cup \ell$$
where $\ell$ is a line.  This yields a subvariety of
$F_{20}$ isomorphic to $\bP^2$;  the lines 
$R\subset \bP^2$ trace out $T_{5,1}$'s on $X_{20}$.  
The corresponding $\rho \in \Pic(F_{20})$ is a 
$(-10)$-class. 
\end{exam}
\begin{exam}[$T_{5,2}$] \label{exam:T52}
There are no quintic scrolls with two ordinary double points in $\bP^5$. 
(This is highly 
unfortunate because a cubic fourfold containing such
a scroll would be rational by Proposition \ref{prop:unirat}.)
Let $\tilde{T}_{5,2}\subset \bP^6$ be the normalization and 
$p\in \bP^6$ a point such that the projection 
of  $\tilde{T}_{5,2}$ from $p$ is $T_{5,2}$. It follows that
$\tilde{T}_{5,2}$ contains four  coplanar points. 
However, these points necessarily lie 
on a conic curve $C\subset\tilde{T}_{5,2}$.
This forces $T_{5,2}$ to be singular along the image of $C$. 
\end{exam}
\begin{exam}[$T_{6,0}$] \label{exam:T60}
This remains to be explored - 
the corresponding $\rho$ 
is {\em not} nodal! The discriminant of the lattice
$\bZ h^2+\bZ T_{6,0}$ is 12.  
\end{exam}
\begin{exam}[$T_{6,1}$] \label{exam:T61}
In this discriminant ($d=18$) the
Fano variety has two square-zero 
classes (bounding the ample cone, by our conjectures). 
In particular, there are no
nodal classes in this case.
\end{exam}
\begin{exam}[$T_{6,2}$] \label{exam:T62}
The scroll $T_{6,2}$ is contained in a general cubic fourfold
$X_{24}$ of discriminant 24. The family of such scrolls in a given
cubic fourfold is parametrized by a K3 surface of degree 6. 
The corresponding class $\rho$ is a nodal $(-2)$-class. 
\end{exam}
\begin{exam}[$T_{6,3}$] \label{exam:T63}
The cubic fourfolds $X_{30}$ 
containing a sextic scroll with three ordinary double points
are codimension 2 in moduli.
The normalization $\tilde{T}_{6,3}\subset \bP^7$ has
6 points lying in a 4-dimensional linear subspace, containing 
the line $\ell $ from which we project.  
These points necessarily are contained in a 
rational normal curve $C
\subset \tilde{T}_{6,3}$
of degree  4.  
The image of $C$ under projection is a quartic plane curve. 
This plane is necessarily contained in $X_{30}$ 
by Bezout's theorem,
so $H^{2,2}(X_{30},\bZ)$ has rank at least 3. 
Let us remark that
the Fano variety $F_{30}$ 
contains a surface isomorphic to $\bP^1\times \bP^1$
and the rulings of the scrolls are given by type $(1,1)$-curves
of this surface (cf. the discussion of $T_{4,1}$ and Example
\ref{exam:F0}).  
\end{exam}
\begin{exam}[Further examples] \label{exam:further}
Essentially the same argument shows 
that there are no scrolls
$T_{6,4}$ (or $T_{6,5}$ or $T_{7,5}$): 
we look at the 8 points on the normalization
$\tilde{T}_{6,4}$ spanning a 5-dimensional 
linear subspace containing
the line of projection $\ell$. 
These points are necessarily contained
on a rational normal curve of degree 5 on $\tilde{T}_{6,4}$.
It projects to a quintic curve in $\bP^3$ with 4 ordinary double
points. This violates Bezout.  
\end{exam}

\section*{Erratum, added March 5, 2010}

\thispagestyle{empty}

In this Erratum we correct two mistakes in the published version of this paper:
\begin{quote}
Rational curves on holomorphic symplectic fourfolds,
{\em Geometric and Functional Analysis} {\bf 11} (2001), no. 6, 1201-1228 
\end{quote}

We are grateful to Antoine Chambert-Loir for pointing out that the
first paragraph of Theorem 4.1 should read as follows:
\begin{quote}
Let $F$ be an irreducible
holomorphic symplectic manifold of dimension $2n$ and $Y$ a
submanifold of dimension $k$.
Assume either that $Y$ is Lagrangian, or that all
of the following hold:
$
\cN_{Y/F}= \Omega^1_Y\oplus \cO^{\oplus 2n-2k}_Y,
$
the restriction of the symplectic form to $Y$ is zero,
and $H^1(\cO_Y)=0$.
Then the deformation space of $Y$ in $F$ is smooth,
of dimension $2n-2k$ {\bf if the last three conditions above hold}.
\end{quote}

We are grateful to Claire Voisin for pointing out that Proposition 7.4
requires an additional hypothesis, and should read as follows:
\begin{quote}
Let $X$ be a cubic fourfold with Fano variety $F$.
Assume that $F$ contains a smooth rational curve $R$
of degree $n$, with corresponding scroll $T_{n,\Delta}$.
Assume that this corresponding scroll $T$ is not a cone
{\bf and has isolated singularities}.
Then there exists a rational map
$$\phi:\bP^4 \dashrightarrow X$$
with
$$
\deg(\phi)
=\binom{n-2}{2}-\Delta
=\frac{(n-2)^2}{4}+\frac{\left(R,R \right)}{2}+1.
$$
\end{quote}
Without this hypothesis the double point computation fails, as is shown
by the example
$$T=\{x^2z+y^2t=0\}.$$

\end{document}